\documentclass[12pt]{article}
\title{Puzzles, Tableaux, and Mosaics}
\author{Kevin Purbhoo  \\
University of British Columbia\\
{\normalsize \tt kevinp@math.ubc.ca}}

\usepackage[mathscr]{eucal}
\usepackage{latexsym}
\usepackage{amssymb}
\usepackage{amsthm}
\usepackage{amsmath}
\usepackage{epsfig}
\usepackage{colordvi}

\newcommand{\CC}{\mathbb{C}}

\newcommand{\NN}{\mathbb{N}}
\newcommand{\RR}{\mathbb{R}}
\newcommand{\ringH}{\mathcal{H}_{d,n}}

\newcommand{\graylozenge}{\begin{picture}(0,0)(0,0)\put(0,0){$\Gray{\blacklozenge}$}\end{picture}\nolinebreak\lozenge}
\newtheorem{theorem}{Theorem}
\newtheorem{lemma}{Lemma}[section]
\newtheorem*{wclemma}{The Wake Crossing Lemma}

\newtheorem{proposition}[lemma]{Proposition}
\newtheorem{corollary}[lemma]{Corollary}

\begin{document}
\maketitle

\begin{abstract}
We define mosaics, which are naturally in bijection with Knutson-Tao
puzzles.  We define an operation on mosaics, which shows they are also
in bijection with Littlewood-Richardson skew-tableaux. Another consequence 
of this construction
is that we obtain bijective proofs of commutativity and associativity
for the ring structures 
defined either
of these objects.  In particular, we obtain a new, easy proof of the 
Littlewood-Richardson rule.
Finally we discuss how our operation is related to other known
constructions, particularly jeu de taquin.
\end{abstract}

\section{Introduction}
\label{sec:intro}

It is well known that the cohomology ring of the Grassmannian 
$\ringH := H^*(Gr(d,n))$ has a natural geometric basis, the 
Schur basis $\{s_\lambda\ |\ \lambda \in \Lambda\}$.  
The structure constants $a_{\nu \mu}^\lambda$ of $\ringH$
in the Schur basis, defined by
$$s_\nu s_\mu = \sum_{\lambda \in \Lambda}a_{\nu \mu}^\lambda s_\lambda\,,$$
are the {\bf Littlewood-Richardson numbers}.  They are non-negative
integers, and also appear
as structure constants in representation ring of $GL(d)$ and the ring of
symmetric polynomials.

Throughout, the integers $0 < d < n$, will be fixed.
The set $\Lambda$ which indexes
the Schur classes $s_\lambda \in \ringH$ can be viewed concretely as either 
the set of $01$-strings of length $n$ with, with exactly $d$ ones,
or as partitions whose diagram fits inside a $d \times n{-}d$ rectangle.
There is a standard way of passing back and forth between these
two: the positions of the zeroes and ones in a $01$-string are respectively
the positions of the horizontal and vertical steps along the boundary
of the corresponding diagram, as shown below. \\
\centerline{\epsfig{file=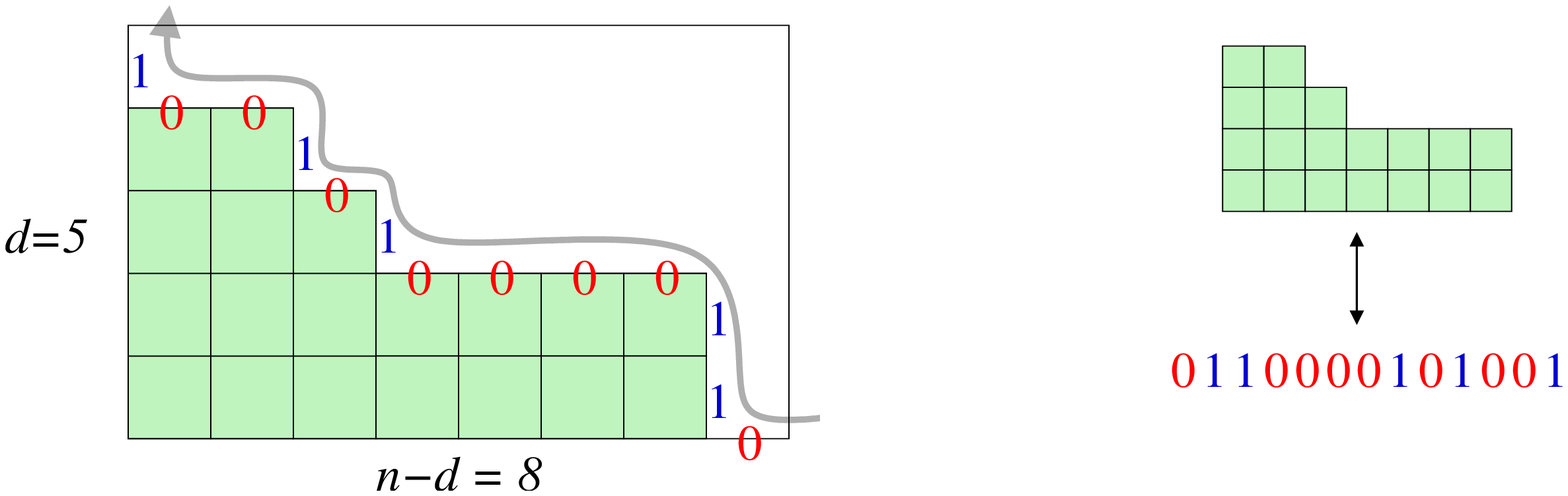, height=1.3in}}

An important problem of the last century has been to find
combinatorial interpretations for the Littlewood-Richardson numbers;
we refer to these collectively as Littlewood-Richardson rules.
In this paper we shall be concerned with two such rules.  
The first, due to Knutson, Tao and Woodward, 
states that $a_{\nu \mu}^\lambda$ can be obtained by 
counting {\em puzzles} with $\nu$, $\mu$ and $\lambda$ on
the boundary \cite{knutson-tao-woodward}.  
The second is the original Littlewood-Richardson rule 
\cite{littlewood-richardson}, which
was proved in \cite{schutzenberger}, and
states that $a_{\nu \mu}^\lambda$ can be obtained by counting 
{\em Littlewood-Richardson tableaux (LR-tableaux)}.  We recall these 
rules and all relevant
definitions in Sections~\ref{sec:puzzles} and~\ref{sec:tableaux}.

The purpose of this paper is to introduce a new 
construction---{\em mosaics}---which will allow us to give new simple
proofs of both of these rules, and provide a bijection between them.
To show that a Littlewood-Richardson rule is correct, one
needs to check two things:
\begin{enumerate}
\item[(i)] that the numbers given by the rule are the
structure constants of a commutative, associative ring;
\item[(ii)] that the Pieri rule holds, that is, multiplication by
special classes which are generators of $\ringH$ behaves correctly.
\end{enumerate}
The Pieri rule is fairly trivial to check for both puzzles and tableaux.
(The rule simply states 
that $a_{(k)\,\mu}^{\ \lambda} = 1$ if the 
diagram $\lambda/\mu$ 
consists of $k$ boxes in $k$ distinct columns, and 
equals zero otherwise, 
where $(k)$ is the partition $k \geq 0 \geq \cdots \geq 0$;
the classes $\{s_{(k)}\}_{k \leq n-d}$ generate $\ringH$.)
Our focus, therefore,  will be on proving commutativity 
and associativity.  N.B. Although we shall prove both here, 
to prove a Littlewood-Richardson rule, associativity alone is
sufficient.

This approach of showing Pieri and associativity was used by 
Knutson, Tao and Woodward to give a proof of 
the {\em hive} formulation of the Littlewood-Richardson rule 
\cite{knutson-tao-woodward:octahedron}.  From this, one can deduce
the puzzle rule.  For LR-tableaux, it is possible to prove commutativity 
and associativity using {\em tableau switching}, as 
defined in \cite{benkart-sottile-stroomer}.  While this approach to
showing commutativity has been discovered many times, the fact that
one can also show associativity (Corollary~\ref{cor:tableauassociate}) 
in this way does not appear to be well documented.
Buch, Kresch, and Tamvakis, also gave simple, parallel proofs 
of the Littlewood-Richardson rule and the puzzle rule along these lines
\cite{buch-kresch-tamvakis}.
Our approach is unique in that the worlds of puzzles and tableaux are 
intertwined, the result of which is that our proofs are surprisingly
short and clean.

Because puzzles have a three-fold symmetry, it will be more
convenient phrase statements of commutativity and associativity
in terms of symmetric Littlewood-Richardson numbers, which
are defined to be
$$a_{\nu \mu \lambda^\vee} := \int_{Gr(d,n)} s_\nu s_\mu s_{\lambda^\vee}\,.$$
These are related to the
ordinary Littlewood-Richardson numbers by the fact that 
$a_{\nu \mu \lambda^\vee} = a_{\nu \mu}^\lambda$,
where $s_{\lambda^\vee}$ is the unique Schur class such that
$$\int_{Gr(d,n)} s_\lambda s_{\lambda^\vee} = 1\,.$$
In terms of partitions, $\lambda^\vee$ is the 
{\em complement} of $\lambda$; in terms of $01$-strings, $\lambda^\vee$ 
is the {\em reverse} of $\lambda$.
Commutativity of $\ringH$ is expressed by the statement that
\begin{equation}
\label{eq:commutativity}
a_{\nu \mu \lambda^\vee} = a_{\mu \nu \lambda^\vee}\,,
\end{equation}
and associativity is expressed as
\begin{equation}
\label{eq:associativity}
\sum_{\kappa \in \Lambda} a_{\nu \mu \kappa^\vee}\, a_{\xi \kappa \lambda^\vee}
= \sum_{\kappa \in \Lambda} a_{\mu \xi \kappa^\vee}\, a_{\kappa \nu \lambda^\vee} \,.
\end{equation}

\subsection{Preliminary examples}
\begin{figure}[htbp]
\begin{center}
\epsfig{file=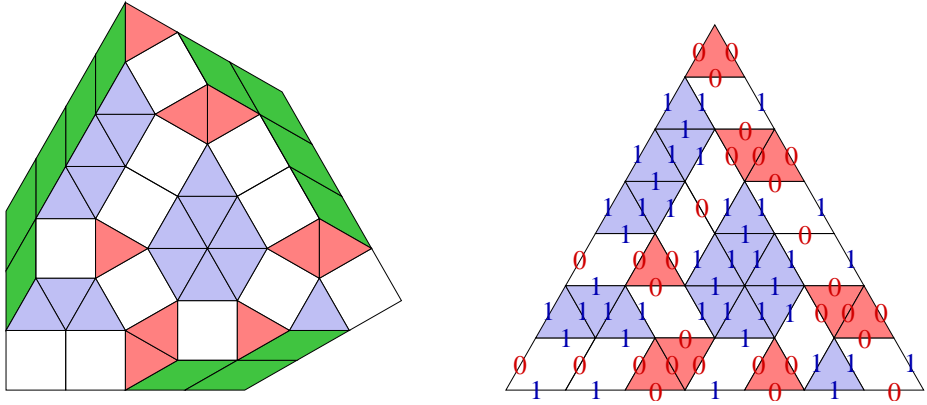, height=2.1in}
\caption{A mosaic with $n=7$, $d=4$, and the corresponding puzzle.}
\label{fig:mosaicex}
\end{center}
\end{figure}

Mosaics are certain tilings
of a hexagonal region.  One can regard a mosaic as a crooked
drawing of a puzzle with some extra rhombi in the corners; in fact, 
there is a straightforward bijection (c.f. Section~\ref{sec:mosaics})
between 
mosaics and puzzles, obtained by removing the extra rhombi, and 
straightening.  
Figure~\ref{fig:mosaicex} illustrates an example of a mosaic and the
corresponding puzzle under this bijection.  

The big advantage of mosaics over puzzles is the fact that 
these extra rhombi
are arranged in the shape of a (slightly distorted) Young diagram.
We define {\em flocks}, which are tableau-like structures on the rhombi,
and an operation called {\em migration} on flocks.  
Our main results state that migration is
an invertible operation which preserves the flock structure.
Thus migration allows us to move flocks around 
inside a mosaic without losing any information.

\begin{figure}[htbp]
\begin{center}
\epsfig{file=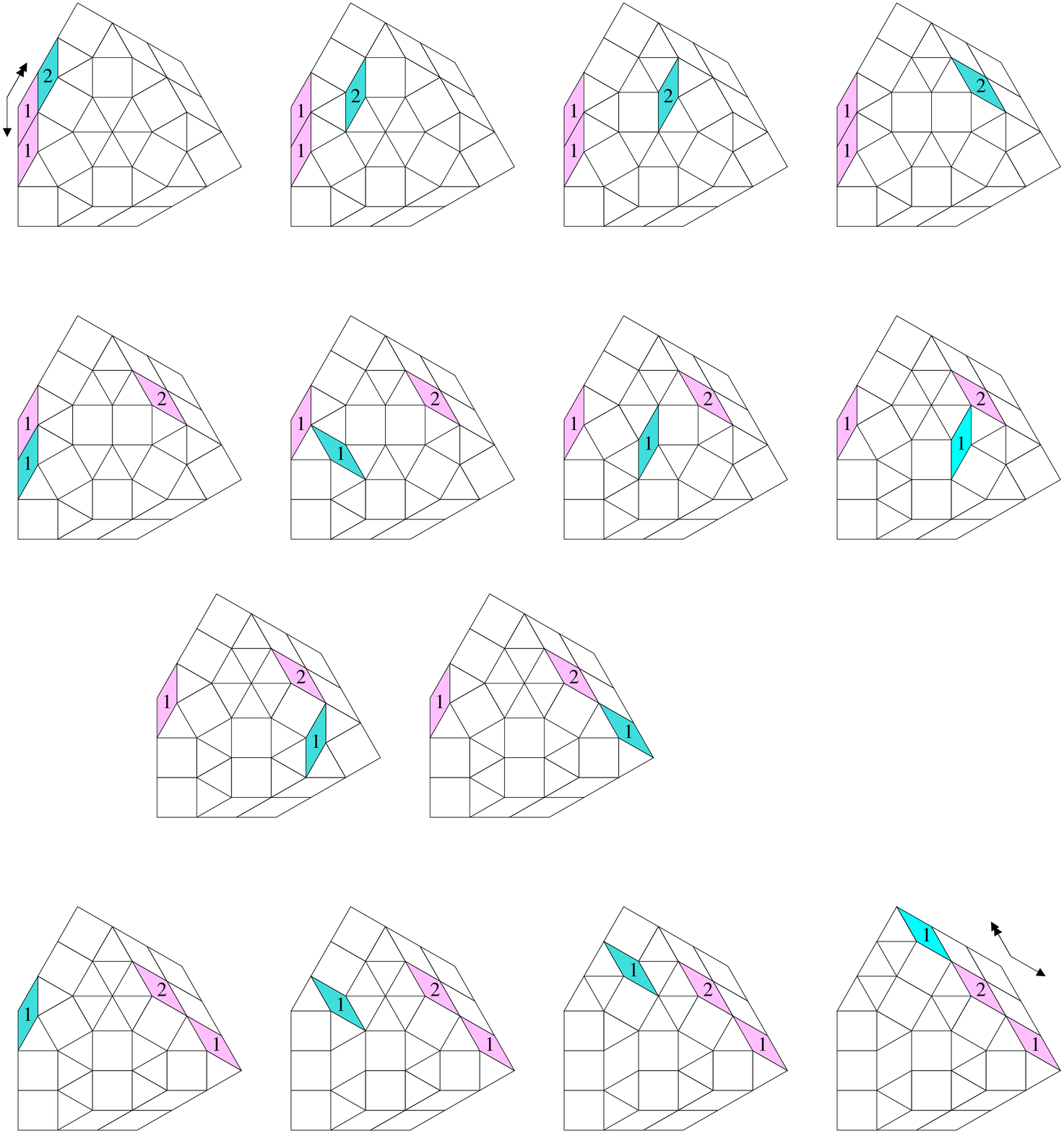, width=6in}
\caption{A flock migrates from the left side of the mosaic to the right.}
\label{fig:migrationex}
\end{center}
\end{figure}

Figure~\ref{fig:migrationex} illustrates an example of a flock on
the left side of a mosaic migrating to the right side.
The result of such a migration is that the flock assumes a shape
which can be interpreted as an LR-tableau, in this case migration
identifies the initial mosaic with the LR-tableau \\
\centerline{\epsfig{file=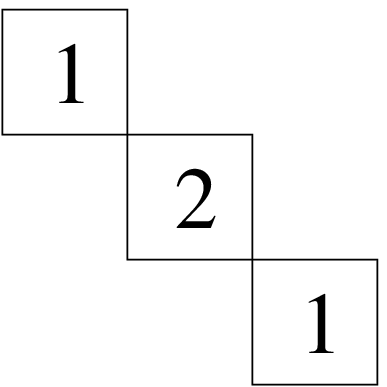, height=.6in}\ \raisebox{5pt}{.}}  
Our main theorems
imply that this process gives a bijection between mosaics 
(or equivalently puzzles) and LR-tableaux. 

One can also perform a sequence of migrations to swap the positions
of two flocks inside a mosaic.  This allows us to deduce that the 
product structure on $\ringH$ defined by puzzles is commutative.
Similar types of arguments allow us to prove commutativity
and associativity for both puzzles and LR-tableaux.  

\subsection{Outline}
This paper is organised as follows.  In Section~\ref{sec:background},
we recall the relevant background information on puzzles and tableaux,
and introduce flocks and mosaics.  The migration operation is defined in 
Section~\ref{sec:operations}.  From here, we give the precise 
statement of our main theorems, and the details of how 
commutativity and associativity follow as a consequence.  We prove
the main theorems in Section~\ref{sec:proofs}.  Finally, in 
Section~\ref{sec:further} we discuss
how migration is related to other known constructions, particularly
to jeu de taquin.  With the exception of this last section, we
have attempted to keep the exposition entirely self-contained.

We would also like to warn the reader that all directions
in this paper are used in a relative, rather than in
an absolute sense.  
When we speak of 
``up'', ``down'', ``left'' and ``right'', these are defined relative 
to a chosen ``right''; hence the reader may sometimes need to rotate
the page so that these are pointing in their usual direction.
Compass directions ``north'', ``east'', ``south'' and
``west'', are used even more loosely.  These are defined relative to 
a chosen basis for $\RR^2$,
which may not be orthogonal or even have the expected orientation.
Hence, east and right are not always 
the same direction, and even if they are, north and up need not be 
the same.  In Sections~\ref{sec:puzzles} and~\ref{sec:tableaux}
these directions have their usual meanings, but thereafter we will
need this additional flexibility.

\subsection{Acknowledgements}
The author is grateful to Allen Knutson for helpful comments on the
paper, and to the Banff International Research Station for providing
the inspirational surroundings in which these ideas were hatched.  
This research was partially supported by NSERC.

\section{Background and definitions}
\label{sec:background}
\subsection{Puzzles}
\label{sec:puzzles}
A {\bf puzzle} of size $n$ is a tiling of an equilateral triangle of 
side length $n$, using the following pieces (each has unit side length)
\begin{center}
\epsfig{file=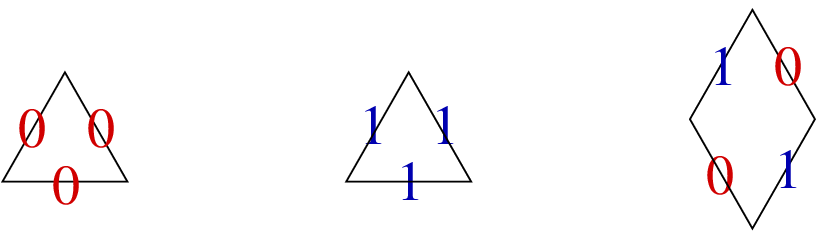,height=0.6in}
\end{center}
subject to the following conditions:
\begin{enumerate}
\item[(i)]
the pieces can be rotated in any orientation, but not reflected;
\item[(ii)]
wherever two pieces share an edge, the numbers on the edge must
agree.
\end{enumerate}

\begin{figure}[htbp]
  \begin{center}
    \epsfig{file=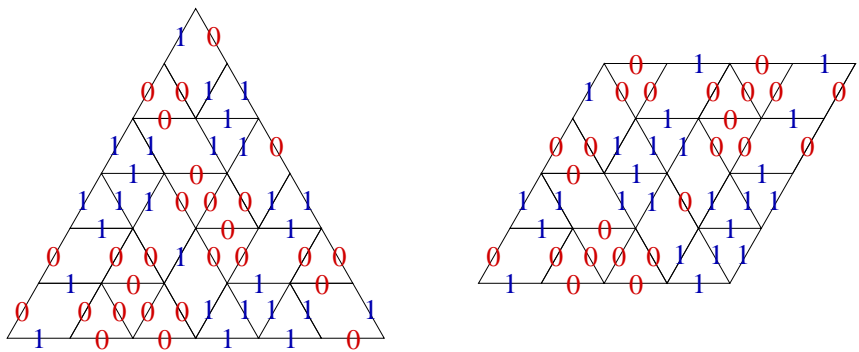,height=1.8in}
    \caption{A puzzle with boundary data $(001101, 010101, 011001)$
and a bipuzzle with boundary data $(0101, 0101, 0011, 1001)$.}
    \label{fig:puzzle}
  \end{center}
\end{figure}

If we read the numbers on the boundary of a puzzle, clockwise starting
from the lower-left corner, we see three strings of $0$s and $1$s,
$\pi$, $\rho$, $\sigma$ of length $n$. See Figure~\ref{fig:puzzle} (left). 
We call $(\pi, \rho, \sigma)$ the {\bf boundary data} of the puzzle.  
For a string $\sigma=\sigma_1\ldots\sigma_n$, let
$\sigma^\vee = \sigma_n \ldots \sigma_1$ denote its reversal.

We recall the following basic fact about puzzles:
\begin{proposition}
If $(\pi, \rho, \sigma)$ is the boundary of a puzzle, then the
strings $\pi$, $\rho$ and $\sigma$ have the same number of $1$s.
\end{proposition}
Henceforth, we shall assume this number equals $d$.

Let $b_{\pi \rho \sigma}$ denote the number of puzzles with
boundary data $(\pi, \rho, \sigma)$.  The puzzle-theoretic
Littlewood-Richardson rule states
$b_{\pi \rho \sigma} = a_{\pi \rho \sigma}$.

A {\bf bipuzzle} is a tiling of the rhombus of side length
$n$ with angles $60^\circ$ and $120^\circ$, satisfying conditions
(i), (ii) and
\begin{enumerate}
\item[(iii)] the number of $1$s appearing on each side of the rhombus 
is equal to the number on any other side.
\end{enumerate}
Again, we assume that this number is equal to $d$.  
The {\bf boundary data} of a bipuzzle is the list of strings
$(\pi, \rho, \sigma, \tau)$ read clockwise from the lower-left
corner.  See Figure~\ref{fig:puzzle} (right).

We will need the following fact.
\begin{proposition}
\label{prop:bisplit}
Every bipuzzle can be split in half into two single puzzles.
\end{proposition}

In other words, there are no rhombi that lie across the line which 
divides the bipuzzle into two equilateral triangles of side length
$n$.  This has an easy ``Green's theorem'' style proof, as in 
\cite[Lemma 5]{knutson-tao:equiv}; we also give an alternate
proof in Section~\ref{sec:proofs}.  Bipuzzles are important for 
associativity (c.f.~\eqref{eq:associativity}), as we have:

\begin{proposition}
\label{prop:countbipuzzles}
The number of bipuzzles with boundary data 
$(\pi, \rho, \sigma, \tau)$ is equal to
$\sum_{\upsilon \in \Lambda} b_{\rho \sigma \upsilon^\vee} \,
b_{\pi \upsilon \tau} =
\sum_{\upsilon \in \Lambda} b_{\rho \sigma \upsilon^\vee} \,
b_{\upsilon \tau \pi}$.
\end{proposition}

\begin{proof}  First note that the two sums are equal,
since by the rotational symmetry of the puzzle rule,
$b_{\pi \upsilon \tau} = b_{\upsilon \tau \pi}$.  
The expression 
$\sum_{\upsilon \in \Lambda} b_{\rho \sigma \upsilon^\vee} \,
b_{\pi \upsilon \tau}$ counts pairs of puzzles
 with boundary data $(\pi, \upsilon, \tau)$ and 
$(\rho, \sigma, \upsilon^\vee)$ for some $\upsilon \in \Lambda$.
But any such pair can be glued along the
$\upsilon$-boundary to form a bipuzzle with boundary
data $(\pi, \rho, \sigma, \tau)$, and by 
Proposition~\ref{prop:bisplit} every bipuzzle arises 
in this way.  
\end{proof}

\subsection{Young tableaux}
\label{sec:tableaux}
A {\bf partition} $\lambda = \lambda_1 \geq \dots \geq \lambda_d \geq 0$ is a
decreasing sequence of integers.  We will assume that all of our partitions
have at most $d$ parts, and that $\lambda_1 \leq n-d$.
The {\bf complement} of $\lambda$ is the partition 
$\lambda^\vee = \lambda^\vee_1 \geq \cdots \geq \lambda^\vee_d \geq 0$, having
parts $\lambda^\vee_k = n-d -\lambda_{d-k}$.

To each partition $\lambda$, we associate
its {\bf diagram}, also denoted $\lambda$.  We adopt the French convention, 
in which the southmost row has $\lambda_1$ boxes, the next higher row 
has $\lambda_2$ boxes, etc.  All rows are assumed to be left justified.  
The {\bf conjugate} of $\lambda$ is the partition $\lambda^t$ whose
diagram has whose diagram has $\lambda_k$ boxes in the $k^\text{th}$ 
column.  A {\bf skew diagram} $\lambda/\mu$ is the diagram consisting of all 
boxes of $\lambda$, which are not also boxes of the subdiagram
$\mu$.  Here we assume,
of course, that $\mu$ and $\lambda$ have the same southwest corner.
We sometimes refer to a partition diagram as a {\bf straight diagram}, 
to emphasize that it is non-skew.


For any skew diagram $\lambda/\mu$, a Littlewood-Richardson tableau 
(or {\bf LR-tableau}),
is a function $f$ which assigns positive integer entries to the boxes 
of $\lambda/\mu$ such that:
\begin{enumerate} 
\item[(i)] the entries are weakly increasing in each row from east to west;
\item[(ii)] 
the entries are strictly increasing in each column from south to north;
\item[(iii)] 
if one forms the word $w_1\ldots w_r$ by listing the entries in
the ordinary reading order (west to east and north to south), then
for all $k \geq 1$, each tail subword $w_s\ldots w_r$ has at least as 
many $k$s as $(k{+}1)$s.
\end{enumerate}

Given an LR-tableau $f$, the {\bf standard order} on the boxes is the 
following: $X <_f Y$ iff either $f(X) < f(Y)$, or $f(X) = f(Y)$ and $X$
is west (or northwest) of $Y$.  This ordering will be key for us.  It
arises in other contexts in tableau theory, for example in
defining Sch\"utzenberger slides.

The {\bf content} of a LR-tableau is the partition $\nu = \nu_1 \geq \dots 
\geq \nu_r$, where $\nu_k$
is the number of $k$s in the tableau.  The {\bf boundary data} of
a LR-tableau of shape $\lambda/\mu$ and content $\nu$ are defined
to be $(\nu, \mu, \lambda^\vee)$.

Let $c_{\nu \mu \lambda^\vee}$ denote the number of LR-tableaux with
boundary data $(\nu, \mu, \lambda^\vee)$.  
The Littlewood-Richardson rule asserts that 
$c_{\nu \mu \lambda^\vee} = a_{\nu \mu \lambda^\vee}$.

An {\bf LR-bitableau} on a skew diagram $\lambda/\mu$
is partition diagram $\kappa$ with $\mu \subset \kappa \subset \lambda$
together with a pair of LR-tableaux $f$ and $g$ on 
$\lambda/\kappa$ and $\kappa/\mu$.  The {\bf boundary data} of this
LR-bitableau are $(\xi, \nu, \mu, \lambda^\vee)$, where $\xi$ and
$\nu$ are the content of $f$ and $g$ respectively.
LR-bitableaux are important for associativity
(c.f.~\eqref{eq:associativity}), since we have:
\begin{proposition}
\label{prop:countbitableaux}
The number of LR-bitableaux with boundary data 
$(\xi, \nu, \mu, \lambda^\vee)$ is 
$\sum_{\kappa \in \Lambda} c_{\nu \mu \kappa^\vee}\, c_{\xi \kappa \lambda^\vee}$.
\end{proposition}

\subsection{Transformed tableaux and flocks}
When we introduce mosaics, we will encounter pictures which look
like skew diagrams that have been subjected to a linear transformation.  
We refer to such a picture as a {\bf transformed diagram}.
We would like to treat transformed diagrams as standard skew diagrams.
However, without knowing the linear transformation, we 
cannot uniquely determine the original skew diagram from the
transformed diagram---there can be (and usually are) up to four 
possibilities.  We therefore define an {\bf orientation} 
on a transformed diagram $\alpha$ to be a pair of vectors 
$E,N \in \RR^2$, such that if 
$\phi$ is the linear transformation given by the 
matrix $(E|N)$, then $\phi^{-1}(\alpha)$ is a standard skew diagram.  
An orientation gives a concrete identification of a transformed diagram
with a standard skew diagram.  To specify and orientation pictorially, 
we draw $E$ with a single headed 
arrow, and $N$ with a double headed arrow, as shown below.
(The arrowheads reflect the fact that the entries of an LR-tableau
weakly increase eastward and strictly increase northward.)
\centerline{\epsfig{file=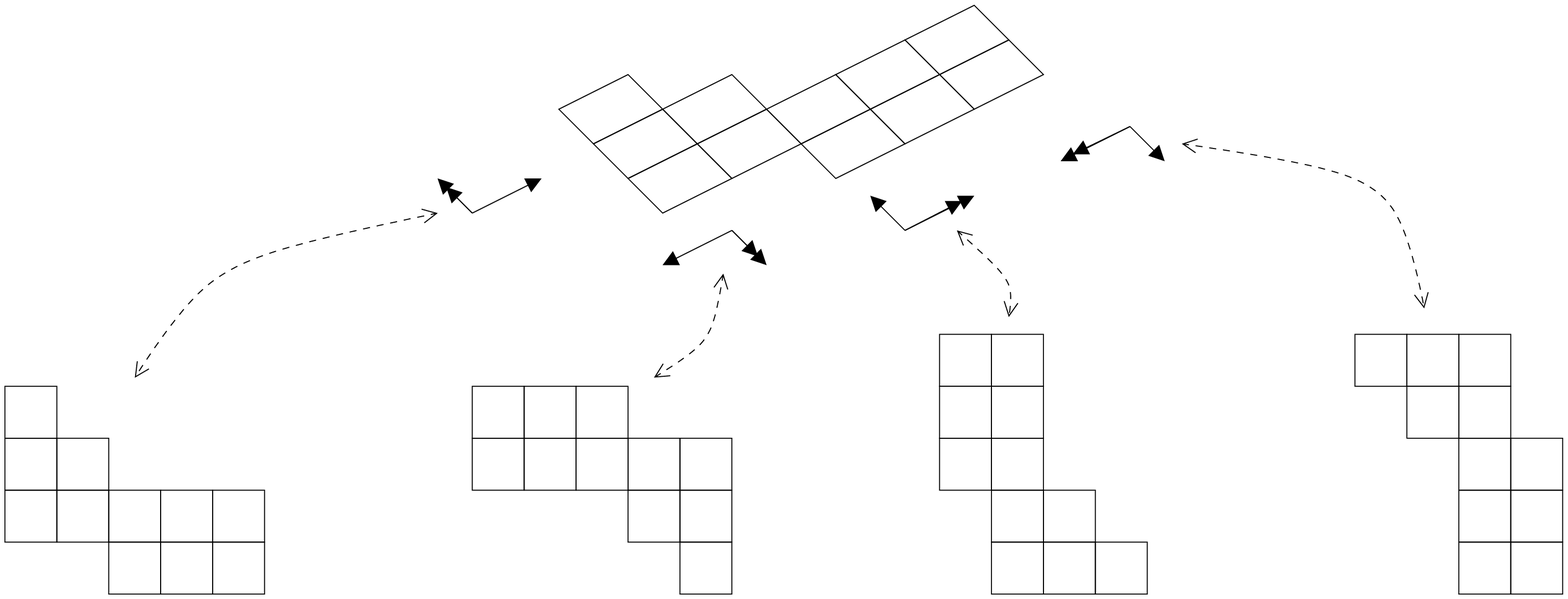,height=1.5in}}
Henceforth we will refer to the directions $E$, $N$, $-E$, $-N$ as east,
north, west and south, respectively.  Once an orientation has been chosen,
we can talk about LR-tableaux, standard order, etc. 
on a transformed diagram---our existing definitions make sense once
the compass directions are thusly redefined.

A {\bf nest} is a designated convex cone in the plane, with 
angle $150^\circ$.
An {\bf orientation} on a nest is a pair of unit vectors $(E,N)$ parallel
to the sides of the cone at an angle of $150^\circ$.  There are
four possible orientations for any nest.
Consider a collection of rhombi of unit side length
with angles $150^\circ$ and $30^\circ$.  We say that the rhombi
are {\bf packed} in the nest, if they lie inside the cone, have
edges parallel to the cone, and are a linear transformation
of a {\em straight} diagram whose corner is the corner of the cone.

\begin{figure}[htbp]
\begin{center}
\epsfig{file=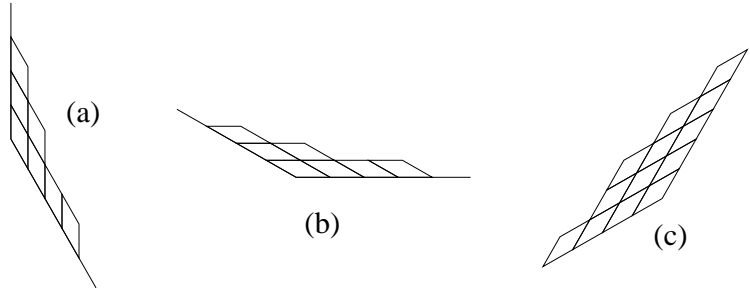, height=1.5in}
\caption{Three collections of rhombi packed in nests; (a) and (b)
are equivalent, while (c) is a complement to both (a) and (b)
when $d=4$, $n=9$.}
\label{fig:equivcomp}
\end{center}
\end{figure}

Given two collections of rhombi $\alpha$ and $\beta$ packed in
two different nests, we say that $\alpha$ is {\bf equivalent}
to $\beta$ if there is an orientation preserving isometry which maps 
one to the other.  We say that $\alpha$ is a {\bf complement} of
$\beta$, written $\alpha = \beta^\vee$,  if there is an orientation 
preserving isometries $\phi_1, \phi_2$ such that $\phi_1(\alpha)$ and 
$\phi_2(\beta)$ 
do not overlap, and together tile the $d \times n{-}d$ parallelogram 
below.  See Figure~\ref{fig:equivcomp}.
\begin{center}
\epsfig{file=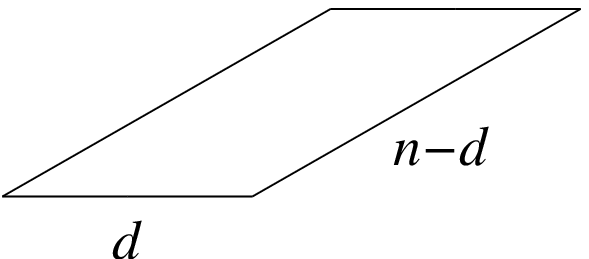, height =.6in}
\end{center}

Given a collection of rhombi $\alpha$ packed in a nest, a triple
$F = (\zeta,(E,N),f)$ is called a {\bf flock} (on $\zeta$) if
\begin{enumerate}
\item[(i)] $\zeta \subset \alpha$ is a subset of these rhombi, in the 
shape of a transformed diagram;
\item[(ii)]
$(E,N)$ is an orientation for both the nest and for $\zeta$;
\item[(iii)] $f:\zeta \to \NN$ is an LR-tableau on $\zeta$ 
in this orientation.
\end{enumerate}
The {\bf content} of a flock is the content of $f$.
If $\alpha$ is the set of all rhombi in the nest, we say that $F$ is 
{\bf accessible} if $\zeta = \alpha/\alpha_1$ for 
some (transformed) straight diagram $\alpha_1 \subset \alpha$.  
A single rhombus $\lozenge \in \alpha$ is accessible if the unique 
flock on $\{\lozenge\}$ is accessible.

Note if the nest has already been given an orientation, we do {\em not}
insist
that a flock have the same orientation.  Thus we can have flocks
of different orientations in the same nest.


\subsection{Mosaics}
\label{sec:mosaics}
\begin{figure}[htbp]
\begin{center}
\epsfig{file=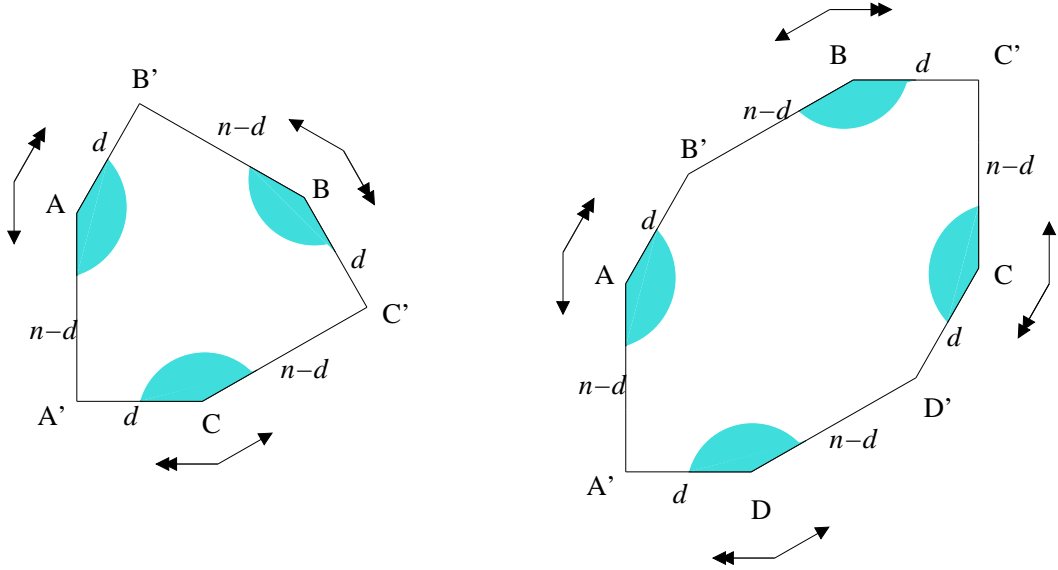, width=5.5in}
\caption{The hexagon $A'AB'BC'C$ and the octagon $A'AB'BC'CD'D$
and their nests, used in mosaics and bimosaics respectively.} 
\label{fig:hexoct}
\end{center}
\end{figure}

Consider the hexagon $A'AB'BC'C$ (vertices listed clockwise) which has 
angles $150^\circ$ at $A$, $B$, $C$,
$90^\circ$ at $A'$, $B'$, $C'$, and side lengths
$A'A = B'B = C'C = d$ and $AB' = BC' = CA' = n-d$.
where $0 < d < n$ are integers.  We designate three nests,
at the corners $A$, $B$ and $C$.  Define unit vectors
$E_A, N_A, E_B, N_B, E_C, N_C$ in the directions of 
$\overrightarrow{AA'}, \overrightarrow{AB'}, 
\overrightarrow{BB'}, \overrightarrow{BC'}, 
\overrightarrow{CC'}, \overrightarrow{CA'}$ respectively, and
fix orientations $(E_A, N_A)$, $(E_B, N_B)$, $(E_C, N_C)$ 
on the nests at $A$, $B$, and $C$ respectively.  See 
Figure~\ref{fig:hexoct} (left).

A {\bf mosaic} is a tiling of this hexagon by the following three shapes:
\\
\parbox{3.5in}{
\begin{enumerate}
\item[(a)] the equilateral triangle with side length 1,
\item[(b)] the square with side length 1,
\item[(c)] the rhombus with side length 1, and 
angles $30^\circ$ and $150^\circ$,
\end{enumerate}
}\hspace{.5in}
\parbox{1.2in}{
\epsfig{file=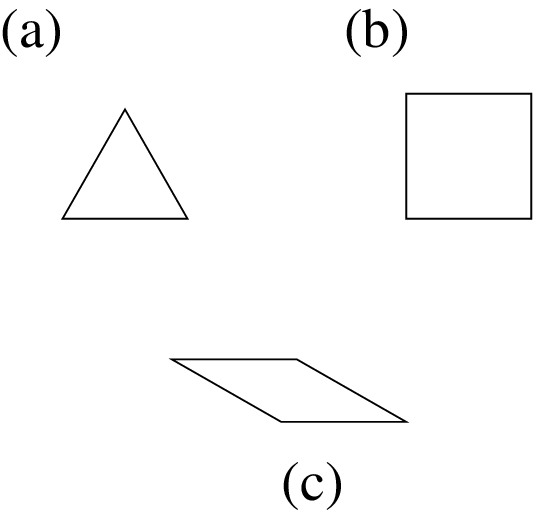, height=1.2in}}
\\
such that all rhombi are packed into the three nests of the hexagon.
The collections of rhombi in nests $A$, $B$ and $C$ are denoted 
$\alpha$, $\beta$, $\gamma$ respectively, and  $(\alpha, \beta, \gamma)$ 
are called the {\bf boundary data} of the mosaic.
The standard orientations on $\alpha, \beta, \gamma$ are the
orientations of the nests $A$, $B$, and $C$ respectively.

There is a natural bijection between mosaics and puzzles, which easiest
to describe physically.  Take a mosaic
and remove all the rhombi.  Now imagine the corners of each square have hinges
so that the angles of the square can flex to become a rhombus.  Grab the 
three corners $A'$, $B'$, $C'$ and pull tight in outward directions.  
The boundary will straighten into an equilateral
triangle of side length $n$, and the squares will become 
$30^\circ$/$60^\circ$-rhombi.  Finally note that some line segments 
will have rotated clockwise, while others will have rotated 
counterclockwise.
This distinction gives the $0$s and $1$s on the puzzle, respectively.
See Figure~\ref{fig:mosaicex}.

The proof that this works comes from an observation of Knutson and 
Tao (unpublished): puzzles 
are actually the projection of a unique piecewise linear surface in 
$\RR^4$, where
the edges point in directions
$(1,0,0,0)$, $(0,1,0,0)$, $(1,-1,0,0)$ if the edge is labelled $0$,
and
$(0,0,1,0)$, $(0,0,0,1)$, $(0,0,1,-1)$ if the edge is labelled $1$.
Mosaics are obtained from a slightly different projection, and the
two projections are related by the operation just described.

This bijection shows us how to compare the boundary data of a puzzle
with the boundary data of a mosaic.  Under this bijection we see that
the shape left by removing $\alpha$ turns into the string of $0$s and 
$1$s corresponding to the walk from $A'$ to $B'$: we get a $0$ for
each unit step west, and a $1$ for each unit step north.
We will say that a mosaic with boundary data 
$(\alpha, \beta, \gamma)$ 
and a puzzle with boundary data $(\pi,\rho, \sigma)$ have the 
{\bf same boundary}, if this correspondence identifies $\alpha$ with $\pi$,
$\beta$ with $\rho$ and $\gamma$ with $\sigma$.
We will write $b_{\alpha \beta \gamma} = b_{\pi \rho \sigma}$
for the number of puzzles with boundary data 
$(\alpha, \beta, \gamma)$.

We can also easily compare the boundary data of an LR-tableau and a 
mosaic.  We will say that an mosaic with boundary data 
$(\alpha, \beta, \gamma)$ and an LR-tableau with boundary
data $(\nu, \mu, \lambda^\vee)$ have the {\bf same boundary}, if the 
standard orientations identify $\alpha$ with $\nu$, $\beta$ with $\mu$
and $\gamma$ with $\lambda^\vee$.  Note that these identifications are
consistent with
the bijection between partitions and $01$-strings described in
Section~\ref{sec:intro}.

We also define bimosaics, which (by the same straightening
operation) naturally correspond to bipuzzles.
Consider the octagon $A'AB'BC'CD'D$ (Figure~\ref{fig:hexoct} (right))
with angles $150^\circ$ at $A$, $B'$, $B$, $C$, $D$, $D'$
and $90^\circ$ at $A'$, $C'$, and side lengths
$A'A = B'B = C'C = D'D = d$, $AB' = BC' = CD' = DA' = n-d$.
We designate corners $A$, $B$, $C$ and $D$ to be nests.  See 
Figure~\ref{fig:hexoct} (right).
A {\bf bimosaic} is a tiling of the octagon by the same
three shapes as a mosaic, with all rhombi packed into the four
nests.  The {\bf boundary data} of a bimosaic is
$(\alpha, \beta, \gamma, \delta)$, the collections of rhombi at 
$A$, $B$, $C$, and $D$ respectively.

\subsection{Canonical constructions}
\begin{figure}[htbp]
\begin{center}
\epsfig{file=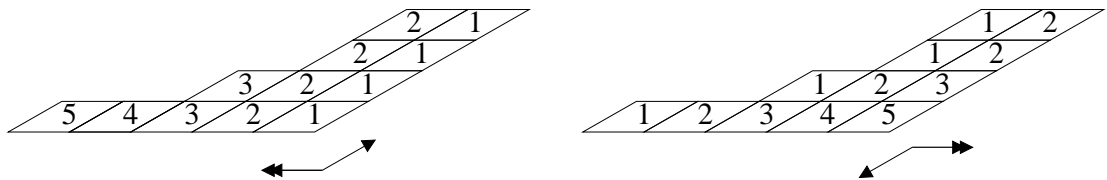, height=.9in}
\caption{The unique flocks on a set of rhombi packed into a nest,
under the orientations $(E,N)$ and $(-E,-N)$.  In both cases
the flock has content $(4,4,2,1,1)$.}
\label{fig:canonicalflock}
\end{center}
\end{figure}

\begin{figure}[htbp]
\begin{center}
\epsfig{file=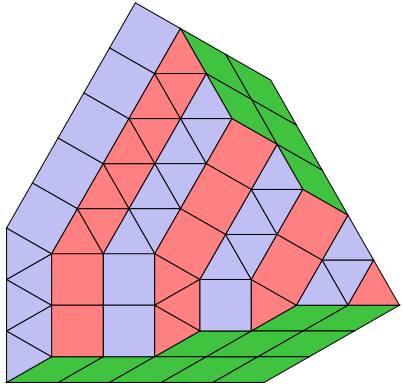, height=2in}
\caption{The unique mosaic when with boundary data 
$(\varnothing, \beta, \beta^\vee)$. Here $n=8$, $d=5$, and 
$\beta = 2 \geq 2 \geq 1 \geq 0 \geq 0$ in the standard orientation.}
\label{fig:canonicalmosaic}
\end{center}
\end{figure}

In a few key situations, certain objects are uniquely constructible.  
The following statements correspond to the fact that multiplying
by identity element in $\ringH$ is trivial.
We omit the proofs as they are straightforward, but illustrative
examples are given in Figures~\ref{fig:canonicalflock} 
and~\ref{fig:canonicalmosaic}.

\begin{lemma}
\label{lem:uniqueflock}
Suppose we are give a nest with an orientation $(E,N)$.
If $\alpha$ is packed in the nest, there is a unique function $f$ 
which makes $F = (\alpha, (E,N), f)$ a flock. 
Moreover for any partition
$\nu$ there is a unique $\alpha$ such that this flock has content
$\nu$.  If the orientation is changed to $(-E,-N)$, the unique flock 
$F^- = (\alpha, (-E,-N), f^-)$ has the same content as $F$.  One
of these two orientations identifies $\alpha$ with the content of $F$.
\end{lemma}

\begin{lemma}
\label{lem:uniquemosaic}
If $\alpha = \varnothing$, then for any $\beta$, there is a unique
mosaic with boundary data $(\alpha, \beta, \gamma)$.  In this
mosaic $\gamma$ is the complement of $\beta$. 
The same is true if the roles of $\alpha, \beta, \gamma$ are
permuted.
\end{lemma}

%

\section{Operations on mosaics}
\label{sec:operations}
\subsection{Migration}
Migration is an operation which makes sense on both mosaics and 
bimosaics.  The operation of migration takes an accessible flock and
moves all of its rhombi in roughly one of four possible compass
directions (e.g. north), so that the flock ends up in a different nest.

The input data for the operation of migration is therefore: (i) a
mosaic or a bimosaic; (ii) an accessible flock 
$F = (\zeta, (E,N), f)$, where $\zeta$ is a collection of rhombi
inside the (bi)mosaic; (iii)
a compass direction, either $\text{north}=N$, 
$\text{east}=E$, $\text{south}=-N$ or $\text{west}=-E$.  The compass
direction
must be consistent with the orientation of the flock, in the following
sense: if the nest in which $F$ is contained opens to the northeast
in this orientation of the flock, the direction must be either north 
or east; otherwise it must be south or west.

In describing migration we will assume that our diagrams are
turned so that the direction of migration is pointing directly
to our right.
In the example of Figure~\ref{fig:migrationex}, the 
direction of migration is north, and thus the page should be turned 
$60^\circ$ clockwise.

For north
or east migration, we begin by locating the rhombus in $\zeta$ which 
is maximal in the standard order.  For south or west migration, we 
locate the minimal one.
Beginning with this rhombus and proceeding in the standard order
(from largest to smallest, or the reverse, whichever is appropriate), 
we perform the following operation with each rhombus $\lozenge$ in 
turn:

Looking to the right of the rhombus $\lozenge$, find the smallest
hexagon in which $\lozenge$ is contained.  (There is always a unique 
minimal hexagon which is {\em entirely} weakly right of the leftmost point 
on $\lozenge$, though sometimes it may seem to be better described as 
above or below $\lozenge$.)
Up to rotation and reflection, 
there are two possible things we could see:
\begin{center}
\epsfig{file=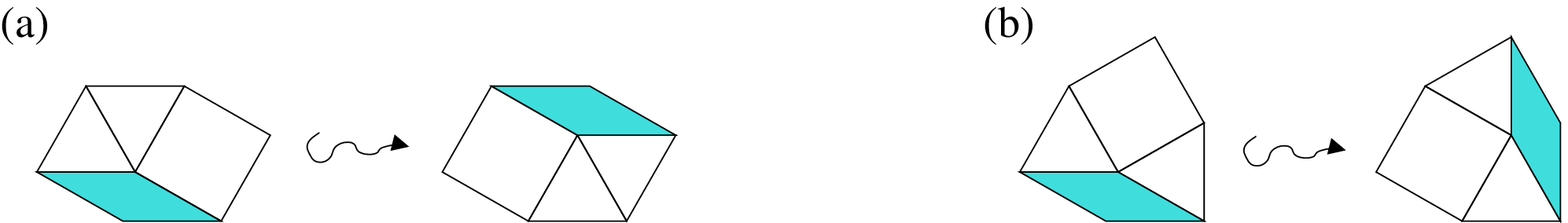, width=4.5in}
\end{center}
If we see (a), the hexagon has a 2-fold rotational symmetry, so we 
rotate the tiling of this hexagon by $180^\circ$.
If we see (b), the hexagon has a 3-fold rotational symmetry, so we 
rotate the tiling 
$120^\circ$, {\em in whichever way causes some edge of the rhombus to end up
exactly horizontal, or exactly vertical}.  (This restriction rules out 
one of two the possible $120^\circ$ rotations.)  In this way, the rhombus
moves roughly to the right (although sometimes it moves strictly vertically).  
Repeat this process, until the rhombus is
packed into a new nest.  Let $\lozenge'$ denote the rhombus in its final 
position, and
let $\zeta' = \{\lozenge'\ |\ \lozenge \in \zeta\}$ denote the collection 
of rhombi from $\zeta$ after migration.

\begin{proposition}
\label{prop:samenest}
All rhombi in $\zeta'$ are in the same nest.
Consequently, $\zeta'$ is a transformed diagram.
\end{proposition}

\begin{proof}
For single (hexagonal) mosaics this is clear, as the rhombi can only ever 
be in 
two of three orientations during the migration process, thereby ruling
out one of the nests.  This argument also works in 
bimosaics for migration towards the opposite nest.  For the other
case (migration towards the near nest), 
Proposition~\ref{prop:bisplit} implies that such a migration 
will behave as it does on a single mosaic.
\end{proof}

In light of this we can specify a direction of migration by specifying
the target nest, rather than the compass direction.  Note that in
a bimosaic, we cannot migrate directly from $A$ to $B$, or 
$C$ to $D$, or vice versa.

There is an induced orientation $(E',N')$ on the transformed
diagram $\zeta'$
which is obtained by rotating the pair $(E,N)$ at most $60^\circ$.
We also have a obvious function $f':\zeta' \to \NN$, which
is defined by $f'(\lozenge') = f(\lozenge)$.  Let $F'$ denote the triple
$(\zeta', (N',E'), f')$,

We can now state our main results.  
\begin{theorem}
\label{thm:isflock}
If $F'=(\zeta', (N',E'), f')$ is the result of migrating a flock 
$F = (\zeta, (N,E), f)$ in some direction inside a mosaic or
a bimosaic, then $F'$ is a flock.  Moreover the map
$\lozenge \mapsto \lozenge'$ preserves the standard order.
\end{theorem}

\begin{theorem}
\label{thm:invertible}
Migration is an invertible operation. 
The inverse operation to
the migration of a flock $F$ from
a nest $A$ to a nest $B$ is the migration of $F'$ from $B$ to $A$.
Inverting, one recovers both the original flock and the 
original mosaic.
\end{theorem}

The proofs are given 
in Section~\ref{sec:proofs}.  The remainder of this section will
be devoted to corollaries of these theorems.

Note: because orientations can become rotated during migration,
migration south is not necessarily the inverse of migration north.
If the orientation is
not rotated (i.e. $E = E'$), then north and south migration are mutually 
inverse, as are east and west migration.
Otherwise, north and west migration are mutually inverse,
as are east and south migration.  In single mosaics, only the latter
case occurs.

\subsection{Bijections between mosaics and tableaux}
\begin{corollary}
\label{cor:mosaictableau}
Migration gives a bijection between mosaics and LR-tableaux with
the same boundary.
\end{corollary}

\begin{proof}
Given a mosaic, with boundary
data $(\alpha, \beta, \gamma)$, give $\alpha$ the orientation
$(-E_A, -N_A)$
and form the unique flock $F = (\alpha, (-E_A,-N_A), f)$.  Let $F$
migrate from $A$ to $B$ (south).
The resulting flock $F'$ is identified with an LR-tableau, which 
(as one can easily check) has the correct boundary data.

The process is reversible.  Given an LR-tableau of shape $\lambda/\mu$,
form the unique mosaic with boundary data 
$(\varnothing, \gamma^\vee, \gamma)$, where
$\gamma^\vee$ is identified with $\lambda$.
Now identify the tableau on $\lambda/\mu$ with a flock $G$
in $\gamma^\vee$.  Let $G$ migrate from $B$ to $A$ (east) to obtain 
a new mosaic.

The two operations are inverse to each other:  migration south is
inverse to migration east, and the reconstruction of the forgotten
data is always unique.
\end{proof}

By giving the flock on $\alpha$ different orientations, we obtain
variants on this bijection.  These are summarised in 
Table~\ref{tbl:variants}.
The second of these bijections is essentially the same as Tao's 
``proof without words'' bijection which appears in \cite{vakil}.  
We discuss this further in Section \ref{sec:jdt}.

\begin{table}[htbp]

\begin{center}
\begin{tabular}{|c|l|}
\hline
Orientation on $\alpha$ & 
Mosaics with boundary data $(\alpha,\beta,\gamma)$ $\leftrightarrow$
\\ \hline \hline
$(-E_A, -N_A)$ & 
LR-tableaux with (same) boundary data $(\nu, \mu, \lambda^\vee)$
\\ \hline
$(E_A, N_A)$ & 
LR-tableaux with boundary data $(\nu, \lambda^\vee, \mu)$
\\ \hline
$(-N_A, -E_A)$ & 
LR-tableaux with  boundary data $(\nu^t, \mu^t, (\lambda^\vee)^t)$
\\ \hline
$(N_A, E_A)$ & LR-tableaux with boundary data $(\nu^t, (\lambda^\vee)^t, \mu^t)$
\\ \hline
\end{tabular}
\end{center}
\caption{Variants on the bijection in the proof 
of Corollary~\ref{cor:mosaictableau}
obtained by changing the orientation on $\alpha$.}
\label{tbl:variants}
\end{table}
%


\subsection{Commutativity}
In the language of puzzles or mosaics, the statement that the 
multiplication on $\ringH$ is commutative, is the following.

\begin{corollary}
\label{cor:mosaiccommute}
There is a bijection between mosaics with boundary data
$(\alpha, \beta, \gamma)$ and those with boundary data
$(\beta, \alpha, \gamma)$.
\end{corollary}

\begin{proof}
Fix any orientations on the nests $A$ and $B$ (not necessarily the
standard ones).  Given a mosaic with boundary data 
$(\alpha, \beta, \gamma)$, form the unique flocks $F$ and
$G$ on $\alpha$ and $\beta$ 
with these orientations.  Migrate
$F$ from $A$ to $C$ to form $F'$.  Migrate $G$ from
$B$ to $A$ to form $G'$.  Finally migrate $F'$ from $C$ to $B$ to
form $F''$.  Since migration does not change the content of
a flock, $F$ and $F''$ have the same content, and thus
$\alpha$ and $\alpha''$ must be equivalent.  Similarly, $\beta$ and 
$\beta'$ are equivalent.  Thus the result of this operation will
have boundary data $(\beta, \alpha, \gamma)$.

The process yields new a orientation on $B$ coming from $F''$,
which is a $120^\circ$ clockwise rotation of the original orientation 
on $A$, and a new orientation on $A$ coming from $F'$, which is
a $60^\circ$ counterclockwise rotation of the original orientation on $B$.
Starting from these orientations, we can reverse the entire process
and recover the original mosaic.  Thus this map is a bijection.
\end{proof}

For tableaux, commutativity of $\ringH$ is expressed as follows.

\begin{corollary}
\label{cor:tableaucommute}
There is a bijection between LR-tableaux with boundary data 
$(\nu, \mu, \lambda^\vee)$ and those with boundary data
$(\mu, \nu, \lambda^\vee)$.
\end{corollary}

\begin{proof}
Given a LR-tableau $f$ on $\lambda/\mu$ with content $\nu$,
form the unique mosaic with boundary data $(\varnothing, \gamma^\vee, \gamma)$,
where $\gamma^\vee$ is identified with $\lambda$.  Let
$F = (\zeta, (E_B, N_B), f)$ be the flock in $\gamma^\vee$ corresponding
to $f$.  Also, let $\eta \subset \gamma^\vee$ be the subset corresponding
to $\nu$, and $G = (\eta, (-E_B, -N_B), g)$ be the unique flock on $\eta$.
Migrate $F$ from $B$ to $A$, then migrate $G$ from $B$ to $A$.
The resulting flocks $F'$ and $G'$ will be such that
$G'$ has boundary data ($\mu, \nu, \lambda^\vee)$.  As the nest
at $B$ is now empty, this process is reversible.
\end{proof}

\subsection{Associativity}

By Proposition~\ref{prop:countbipuzzles},
associativity takes the following form for mosaics.
\begin{corollary}
There is a bijection between bimosaics with boundary data
$(\alpha, \beta, \gamma, \delta)$ and those
with boundary data
$(\delta, \alpha, \beta, \gamma)$.
\end{corollary}

\begin{proof}
There are many ways to do this.  One simply has to form flocks
on $\alpha$, $\beta$, $\gamma$, and $\delta$ and shuffle them
around appropriately inside a bimosaic.   Again, the invertibility 
of such a shuffle operation relies on the fact that orientations rotate in 
predictable way under migration.
\end{proof}

The situation for tableaux is a little less clean, since we
do not have the manifest rotational symmetry of the puzzle rule.
The following assertion for LR-bitableaux does not precisely
show associativity: by Proposition~\ref{prop:countbitableaux}
it translates to the statement that
$x(yz) = z(xy)$.  This plus commutativity gives associativity.

\begin{corollary}
\label{cor:tableauassociate}
There is a bijection between LR-bitableaux with boundary
data $(\xi, \nu, \mu, \lambda^\vee)$ and those with boundary data
$(\nu, \mu, \xi, \lambda^\vee)$
\end{corollary}

\begin{proof}
An LR-bitableau with boundary data $(\xi, \nu, \mu, \lambda^\vee)$
is a pair of LR-tableaux $f$ on $\lambda/\kappa$ with content $\xi$
and $g$ on $\kappa/\mu$ with content $\nu$.
Given a such an LR-bitableau, form the unique mosaic with boundary 
data $(\varnothing, \gamma^\vee, \gamma)$, where $\gamma^\vee$ is identified 
with $\lambda$.  
Let
$F = (\zeta, (E_B, N_B), f)$ be the flock in $\gamma^\vee$ corresponding
to $f$. Let  $G = (\eta, (E_B, N_B), g)$ be the flock corresponding
to $g$. Let $\theta \subset \gamma^\vee$ correspond to $\mu$ and let 
$H = (\theta, (-E_B, -N_B), h)$ be the unique flock on $\theta$.

First migrate $F$ from $B$ to $A$. Then migrate $G$ from $B$ 
to $C$ to form $G'$. 
Next, migrate
$H$ from $B$ to $A$ to produce $H'$, clearing out the nest at $B$.
Finally migrate $G'$ from $C$ to $A$ to produce $G''$ 
The resulting flocks $G''$ and $H'$ together determine an
LR-bitableau with boundary data
$(\mu, \nu, \xi, \lambda^\vee)$.  As in the proof of commutativity,
this is a bijection.
\end{proof}

\section{Proof of main theorems}
\label{sec:proofs}

Once again, we describe our operations under the assumption that
that the direction of migration is to our right.

\begin{figure}[htbp]
\begin{center}
\epsfig{file=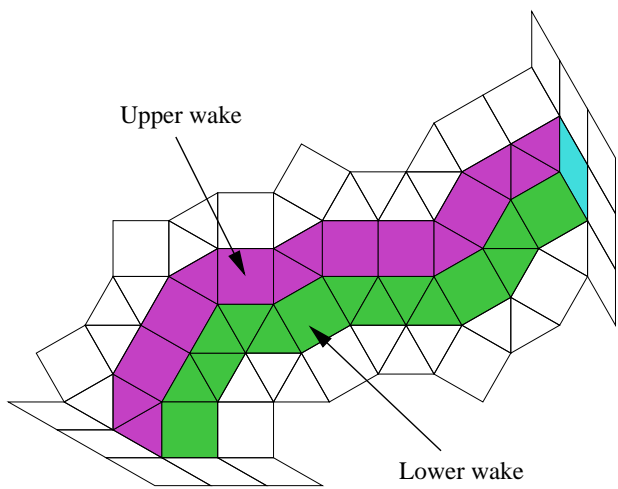, height=2.5in}
\caption{The upper and lower wakes, after the highlighted rhombus
$\lozenge$ has migrated to the right.}
\label{fig:wake}
\end{center}
\end{figure}

Consider the journey of a single rhombus $\lozenge$ from one nest to 
another in the process of migration.  The set of 
tiles (triangles and squares) that are displaced by this journey
form a path from the initial position of $\lozenge$ to 
the final position $\lozenge'$.  We call the set of all these
tiles the {\bf wake} of
$\lozenge$.  This path is everywhere two tiles wide, and 
can be subdivided into two parallel paths, one directly above the 
other, each one tile wide.  We call these the {\bf upper wake} and
{\bf lower wake} of $\lozenge$.  The line which divides them
is called the {\bf midline} of $\lozenge$.  See Figure~\ref{fig:wake}.

Later in the process, when another rhombus $\graylozenge$ takes
its journey it may or may not displace some of the tiles in the 
wake of $\lozenge$.  We will say that $\graylozenge$
{\bf does not disturb} the upper (resp. lower) wake of $\lozenge$
if the journey of $\graylozenge$ does not cause any tiles from
the upper (resp. lower) wake of $\lozenge$ to move.  
We say that the upper (resp. lower) wake of $\lozenge$ is
{\bf intact} for $\graylozenge$, if every rhombus which journeys
after $\lozenge$ and before $\graylozenge$ does not disturb the
upper (resp. lower) wake of $\lozenge$.  

\begin{wclemma}
Suppose $\lozenge$ travels before $\graylozenge$ during migration.
If the upper (resp. lower) wake of $\lozenge$ is intact for 
$\graylozenge$, and $\graylozenge$ begins above (resp. below) the midline 
of $\lozenge$, then the journey of $\graylozenge$ does not cross
the midline.  In particular if the lower (resp. upper) wake
of $\lozenge$ is also intact for $\graylozenge$, then $\graylozenge$
does not disturb the lower (resp. upper) wake of $\lozenge$.
\end{wclemma}

\begin{proof}
There are only a few ways that 
that $\graylozenge$ can plausibly approach the upper wake 
of $\lozenge$.
\begin{center}
\epsfig{file=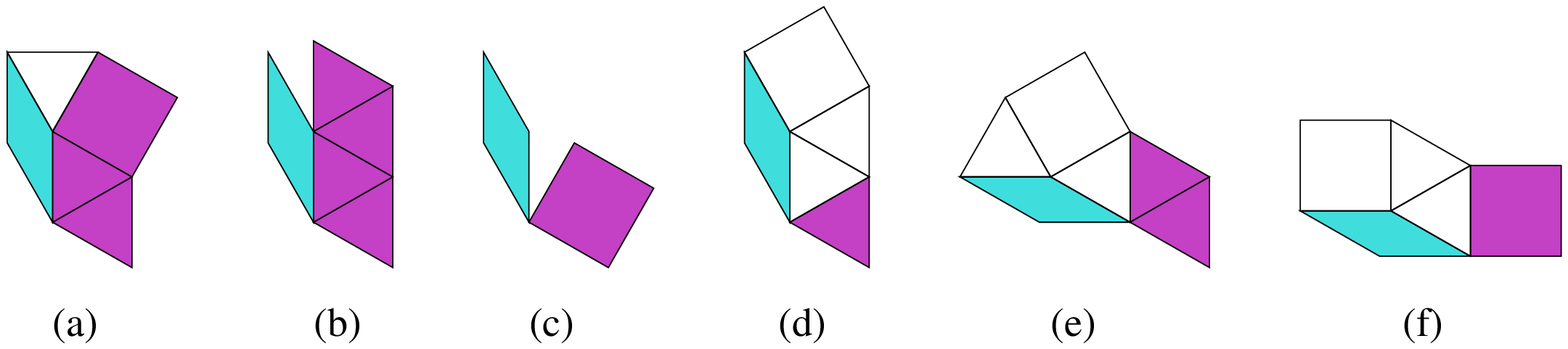, height=1in}
\end{center}
Of these (a) is the only case where $\graylozenge$ disturbs the
upper wake.  The situations in (b) and (c) are impossible, since
the tiling cannot be completed, and in (d), (e) or (f) (or anything
similar), the wake is not immediately disturbed.
One step after (a) occurs, $\graylozenge$ has one edge perpendicular
to the midline, and subsequently retains this property for as long as 
it touches the midline.  As a result $\graylozenge$
can never cross the midline. By symmetry the same is true for the lower wake.
\end{proof}

The upshot is that the rhombi can't get too badly out
of order during migration.  Using the Wake Crossing Lemma,
we deduce Theorem~\ref{thm:isflock}.

\begin{proof}[Proof of Theorem \ref{thm:isflock}]
We will assume that the direction of migration is east or north
and that $E$ is clockwise of $N$.
The proof is essentially the same for south or west migration,
with the standard order reversed.  If $E$ is counterclockwise of $N$,
we swap the roles the upper and lower wake.

We must show that $F'$ satisfies conditions (i), (ii) and (iii)
in the definition of an LR-tableau.
Note that (i) is immediate,
simply because migration proceeds in the standard order.

Let $\lozenge_1 >_f \lozenge_2 >_f \dotsb >_f \lozenge_s$
be the rhombi in $\zeta$ for which $f(\lozenge_i) = k$,
and let $\graylozenge_1 >_f \graylozenge_2 >_f  
\dotsb >_f \graylozenge_t$  ($t \geq s$) be those for 
which $f(\graylozenge_j) = k-1$.

Note that $\lozenge_{i+1}$ is above $\lozenge_i$ in $\zeta$
and the upper wake of $\lozenge_i$ is intact for $\lozenge_{i+1}$.
Therefore, by the Wake Crossing Lemma, $\lozenge'_{i+1}$ is above 
$\lozenge'_i$ in $\zeta'$.  
This implies that $F'$ satisfies condition (ii), and that the map
$\lozenge \mapsto \lozenge'$ preserves standard order.  It also implies that 
the lower wake of each $\lozenge_i$ is intact for $\lozenge_{i+2}, \dotsc, 
\lozenge_s, \graylozenge_1$.  Now, $\graylozenge_i$ is below
$\lozenge_i$ in $\zeta$.  Thus, by the Wake Crossing Lemma
 and a simple induction we deduce
that for each $i$ the lower wake of $\lozenge_i$ is intact for 
$\graylozenge_i$, and thus $\graylozenge'_i$ is below $\lozenge'_i$
in $\zeta'$.  This implies condition (iii).
\end{proof}

\begin{proof}[Proof of Theorem~\ref{thm:invertible}]
It is clear that the journey of each rhombus is invertible, and
the inverse journey is just migration to the original nest.  Since
migration preserves the standard order, performing the inverse
migration will send the rhombi back in the correct order.
\end{proof}

\begin{proof}[Proof of Proposition~\ref{prop:bisplit}]
Under the bijection between bipuzzles and bimosaics, the
line which divides a bipuzzle into two equilateral triangles
corresponds to an edge-path from
$B'$ to $D'$ consisting only of steps perpendicular to $AB'$ or
$B'B$.  Call this path the dividing line of the bimosaic.
The proposition is therefore equivalent to the assertion
that every bimosaic has a dividing line.

Suppose that we have a bimosaic with a dividing line. 
Consider what happens when a single rhombus migrates to the opposite
nest (from $A$ to $C$, or $B$ to $D$, or vice-versa).  
It is easy to see that the when the rhombus
crosses the dividing line, the dividing line simply changes at the
point where the rhombus crosses.  Thus we deduce that the existence of 
a dividing line is invariant under migration to the opposite nest.

Given any bimosaic with boundary data $(\alpha, \beta, \gamma, \delta)$,
form flocks on $\beta$ and $\gamma$ and migrate these to the opposite
nests (in either order).  The result will be a bimosaic with
boundary data $(\alpha', \varnothing, \varnothing, \delta')$, and it is easy
to check that every such bimosaic has a dividing line.  We deduce that
the original bimosaic has a dividing line, as required.

Note that although Proposition~\ref{prop:samenest} 
is implicitly used here, the argument is non-circular as the case
in question (i.e. migration to the opposite nest)
does not rely on Proposition~\ref{prop:bisplit}.
\end{proof}

\section{Further connections}
\label{sec:further}
\subsection{Migration and jeu de taquin}
\label{sec:jdt}

In this section, we will see how migration in mosaics corresponds to
jeu de taquin in tableaux.  We assume familiarity with jeu de taquin
and Sch\"utzenberger slides \cite{schutzenberger} and refer the reader 
to \cite{fulton:tableaux}.  Further algorithms can be found in
\cite{benkart-sottile-stroomer} and the references therein.
We will omit proofs here, since they are long but relatively 
straightforward.

\begin{figure}[htbp]
\begin{center}
\epsfig{file=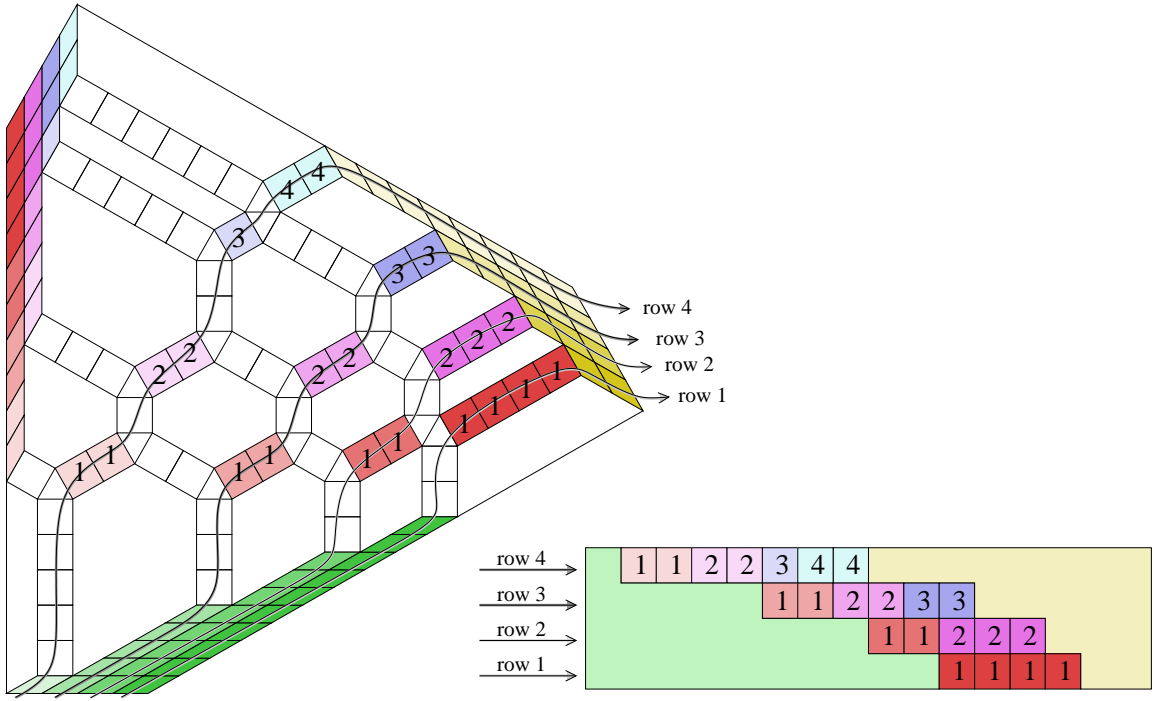, width=6in}
\caption{Bijection between mosaics and Littlewood-Richardson tableaux.}
\label{fig:tao}
\end{center}
\end{figure}

To establish such a connection we need a more explicit description
of the bijection between puzzles and tableaux.  Figure \ref{fig:tao}
is Tao's ``proof without words'' bijection which has been reformulated
in terms of mosaics and Littlewood-Richardson tableaux.  The original
picture, which appears in \cite{vakil}, uses puzzles and an alternate
formulation of the Littlewood-Richardson 
rule \cite[Cor. 5.1.2]{fulton:tableaux}.

\begin{proposition}
\label{prop:samebijection}
The following bijections between mosaics and Littlewood-Richardson
tableaux coincide:
\begin{enumerate}
\item[(i)] Tao's bijection as formulated in Figure~\ref{fig:tao};
\item[(ii)] the bijection obtained by $\alpha$ orientation 
$(E_A, N_A)$ and migrating from $A$ to $B$ (second bijection in 
Table~\ref{tbl:variants});
\item[(iii)] the bijection obtained by giving $\alpha$ orientation 
$(-E_A, -N_A)$ and migrating from $A$ to $C$.
\end{enumerate}
\end{proposition}

This correspondence allows us to see the connection with jeu de taquin.

\begin{proposition}
\label{prop:slides}
Under the bijection of Proposition \ref{prop:samebijection} the migration
of a single rhombus from $B$ to $C$ (or vice-versa) corresponds to
a Sch\"utzenberger slide.
\end{proposition} 

\begin{corollary}
\label{cor:rotate}
Under the bijection of Proposition \ref{prop:samebijection},
rotating a mosaic $120^\circ$ counterclockwise corresponds to the
transformation $f \mapsto \tilde g$ of LR-tableaux described below.
\begin{enumerate}
\item Begin with an LR-tableau $f$ with boundary data $(\nu, \lambda^\vee, \mu)$.
\item Replace each entry by its position in the reverse of the
standard order and rotate this by $180^\circ$.  The result will have
shape $\mu^\vee/\lambda^\vee$.  Call this tableau $\tilde f$.
\item Let $g$ be the unique straight LR-tableaux of shape $\lambda^\vee$.
\item Perform reverse jeu de taquin on $g$: slide the boxes of $\tilde f$, in 
order, through $g$, to produce a new tableau $\tilde g$ with boundary data 
$(\lambda^\vee, \mu, \nu)$.
\end{enumerate}
\end{corollary}

\begin{proof}
By Proposition~\ref{prop:slides} the tableau $\tilde g$ is the result of 
migrating rhombi from $B$ to $C$ with the orientation $(E_B, N_B)$.
\end{proof}

Thus we see that any migration corresponds to some sequence 
of jeu de taquin operations, using Proposition~\ref{prop:slides}
and Corollary~\ref{cor:rotate}.  Hence, questions about migration can 
be translated into questions about jeu de taquin.  We will not attempt
to fully explore the consequences of this correspondence here, but since
so much is already known about jeu de taquin, this is a powerful fact.
For example, the following statement is not at all obvious from the
results in Section~\ref{sec:operations}, but can be shown using jeu
de taquin methods.

\begin{proposition}
\label{prop:involutions}
Each of the bijections described in the proofs of 
Corollaries~\ref{cor:mosaiccommute} and~\ref{cor:tableaucommute} 
is an involution.
\end{proposition}

\subsection{Open questions}
\label{sec:questions}
The reader may at first wonder whether the map 
$\lozenge \mapsto \lozenge'$
is a picture in the sense of Zelevinsky.  It is easily seen that this
is not the case; however, 
since $(f,f')$ is a pair of LR-tableaux of the same content, this
pair corresponds to a picture,
by Zelevinsky's generalisation of the Robinson-Schensted-Knuth
correspondence \cite{zelevinsky}.  Is there
an easy way to describe this picture in terms of migration.

One aspect of this picture which is a bit perplexing is the number
of choices available in making any of the constructions which
prove commutativity and associativity.  Even after one chooses orientations 
on the flocks, there is more than one way to move the flocks around
inside a (bi)mosaic so that they end up in the correct final positions.
Each of these different choices could conceivably give rise to
different commutors and associators (i.e. bijections which prove
commutativity and associativity).  However, the fact that questions
about these operations can be reformulated in terms of jeu de taquin 
limits the number of possibilities.  It would be nice to be able to
see that certain bijections were independent of choices without
resorting to the fact that this is true for jeu de taquin.


One can consider the groupoid whose objects are 
(bi)mosaics with flocks, and whose arrows are sequences of migrations.
What is this groupoid?  Is there monodromy, i.e. when the flocks return 
to their original positions in their original orientations, do we
necessarily get back to the same mosaic?  It should be possible to
answer this question using the correspondence with jeu de taquin,
but again, it would be preferable to somehow get an answer 
directly.

Finally, puzzles, and therefore mosaics, have a geometric 
interpretation~\cite{vakil} in terms of degenerations of Richardson 
varieties.  Can migration be understood in terms of this picture?
If not, perhaps migration can be geometrically interpreted in
a different way.  If this is possible, it may be the nicest way
to understand why certain operations are independent of choices,
as well as answer the monodromy question.

\end{document}